\documentclass[12pt,handout]{article}%
\usepackage[paperwidth=8.5in,left=0.75in,right=0.75in,top=1.0in,bottom=1.0in, paperheight=11.0in]{geometry}
\usepackage{tikz} 
\usetikzlibrary{positioning}

\usepackage{graphicx}
\usepackage{authblk}
\usepackage{amssymb}
\usepackage{epstopdf}
\usepackage{subfigure}  
\usepackage[sans]{dsfont}
\usepackage[applemac]{inputenc}
\usepackage[english]{babel}
\usepackage{latexsym}
\usepackage{mathrsfs}
\usepackage{amscd}
\usepackage{graphicx}
\usepackage{color}
\usepackage{float}
\usepackage{amsmath}
\usepackage{amsfonts}
\numberwithin{equation}{section}
\usepackage{enumerate}
\usepackage{amsthm}
\usepackage{algorithm}
\usepackage{algorithmic}
\usepackage[numbers,sort]{natbib}
\usepackage[bookmarks=true,colorlinks=true,linkcolor={blue},urlcolor={blue}, citecolor={blue},pdfstartview={XYZ null null 1.22}]{hyperref}%

\providecommand{\keywords}[1]{\textbf{{Key words.}} #1}

\title{A data augmentation strategy for deep neural networks with application to epidemic modelling}
 \date{}
\author{Muhammad Awais\footnote{Department of Mathematics and Computer Science, University of Ferrara, Italy, e-mail: wsammm@unife.it}, Abu Safyan Ali\footnote{Department of Mathematics and Computer Science, University of Ferrara, Italy, e-mail: abusafyan.ali@unife.it}, Giacomo Dimarco\footnote{Department of Mathematics and Computer Science, University of Ferrara, Italy, e-mail: giacomo.dimarco@unife.it}, Federica Ferrarese\footnote{Department of Mathematics and Computer Science, University of Ferrara, Italy, e-mail: federica.ferrarese@unife.it}, Lorenzo Pareschi\footnote{Department of Mathematics and Computer Science, University of Ferrara, Italy, and Maxwell Institute for Mathematical Sciences and Department of Mathematics Heriot-Watt University, Edinburgh, UK, e-mail: lorenzo.pareschi@unife.it},}

\begin{document}
 
	\maketitle
\abstract{
	In this work, we integrate the predictive capabilities of compartmental disease dynamics models with machine learning’s ability to analyze complex, high-dimensional data and uncover patterns that conventional models may overlook. Specifically, we present a proof of concept demonstrating the application of data-driven methods and deep neural networks to a recently introduced Susceptible-Infected-Recovered type model  with social features, including a saturated incidence rate, to improve epidemic prediction and forecasting. Our results show that a robust data augmentation strategy trough suitable data-driven models can improve the reliability of Feed-Forward Neural Networks and Nonlinear Autoregressive Networks, providing a complementary strategy to Physics-Informed Neural Networks, particularly in settings where data augmentation from mechanistic models can enhance learning.  This approach enhances the ability to handle nonlinear dynamics and offers scalable, data-driven solutions for epidemic forecasting, prioritizing predictive accuracy over the constraints of physics-based models. Numerical simulations of the  lockdown and post-lockdown phase of the COVID-19 epidemic in Italy and Spain validate our methodology.}

\keywords{Data Driven models in epidemiology, Deep learning, Data augmentation, Feed-Forward Neural Networks, Nonlinear Autoregressive Networks, COVID 19 data}

\maketitle

\section{Introduction}\label{sec1}

The COVID-19 pandemic has underscored the need for accurate and computationally efficient predictive models that not only describe the spread of infectious diseases but also provide fast and reliable forecasts to guide public health interventions. A key challenge in epidemic modeling is balancing mechanistic realism with computational tractability. Among the most widely used approaches are compartmental models, such as the classical Susceptible-Infected-Recovered (SIR) models, which describe disease dynamics by partitioning the population into epidemiologically relevant compartments^^>\cite{capasso1978,hethcote2000,kermack1927}. Specifically, these models categorize individuals as susceptible (capable of contracting the disease), infected (currently infected and capable of transmitting the disease), and recovered (immune or isolated from the population, i.e. people that no longer contribute to disease transmission). Despite their simplicity and effectiveness in capturing some fundamental epidemic trends, classical SIR models assume constant transmission and recovery rates, limiting their ability to represent the complex, nonlinear, and time-varying nature of real-world epidemics.

To address these limitations, extended versions of the SIR model have been proposed in the recent past, incorporating time- and state-dependent transmission rates, saturated incidence functions, and heterogeneous population behaviors^^>\cite{albi2022kinetic,brauer2019,diekmann2000}. For example, models incorporating saturated or nonlinear incidence rates^^>\cite{dimarco2021kinetic,Vespignani,Diekmann}, as well as kinetic approaches to contact dynamics^^>\cite{bertaglia2021spatial}, have proven effective in capturing the impact of social distancing measures during the pandemic. These refinements enhance realism and allow for better representation of the diverse factors influencing epidemic trajectories.

In parallel, machine learning techniques, particularly neural networks, have emerged as powerful tools for studying epidemic spread^^>\cite{brunton2022data,baker2018brochure,la2020epidemiological}. By leveraging their ability to approximate complex, nonlinear relationships, neural networks complement traditional modeling approaches and extend their applicability. Recent advancements, such as Physics-Informed Neural Networks (PINNs), integrate mechanistic knowledge of underlying physical or biological processes directly into the training phase of neural networks^^>\cite{DeRyck_Mishra_2024}. This integration ensures model robustness and consistency with known dynamics, even when data are scarce or uncertain. For instance, PINNs have been successfully applied to solve forward and inverse problems in epidemic modeling, ensuring that predictions adhere to underlying physical laws^^>\cite{bertaglia2022asymptotic,jin2023asymptotic,ouyoussef2024physics, berkhahn2022physics, han2024approaching}.

Feed-Forward Neural Networks (FNNs) and Nonlinear Autoregressive Networks (NARs) also offer promising approaches for epidemiological predictions^^>\cite{Goodfellow-et-al-2016,amendolara2023lstm,li2020recurrent}. FNNs map inputs directly to outputs through layers of nonlinear transformations, making them effective for tasks such as interpolation and short-term forecasting. This approach has been applied across various domains, including early recolonization detection in Parkinson’s disease, biological parameter estimation in disease models, COVID-19 data analysis, musculoskeletal disorder risk assessment, and statistical modeling for selecting optimal FNN architectures^^>\cite{fadavi2024early, rao2010estimation, mcinerney2024statistical, pashaei2021training, chen2000new, liu2024combined}. However, FNNs do not account for temporal dependencies, limiting their ability to model the evolution of epidemics over time.

Conversely, NARs are specifically designed to process sequential data, making them well-suited for time-series forecasting in epidemic modeling. By incorporating previous values of the time series, NARs retain information about past states, enabling them to capture temporal dependencies in epidemic dynamics.  These networks have been successfully applied in the recent past to forecast infection peaks, estimate recovery rates, and simulate intervention scenarios^^>\cite{li2020recurrent, liao2021sirvd, nabi2021forecasting, jamshidi2020artificial, puleio2021recurrent, fokas2020mathematical,bontempi2013machine}.

In this work, we introduce a hybrid epidemic forecasting approach that integrates mechanistic epidemiological models with data-driven FNNs and NARs, creating a surrogate model capable of generalizing across different epidemic scenarios. This approach enables fast, data-driven epidemic forecasting while preserving the interpretability of traditional mechanistic models. The resulting methodology is applied to analyze the lockdown and post-lockdown phase of the COVID-19 epidemic in Italy and Spain, leveraging real-world data and enhanced computational strategies to improve predictive accuracy. This framework enables us to augment real-world epidemic data with synthetic simulations generated from a generalized SIR model, which includes a saturated incidence rate to better reflect intervention effects and behavioral changes. The main idea consists in addressing the limitations posed by sparse real-world data by generating additional synthetic information using a generalized SIR data driven model. Specifically, we focus on an extension known as the social-SIR model with a saturated incidence rate, referred to here as the generalized SIR model in the following. This model modifies the transmission dynamics^^>\cite{dimarco2021kinetic} by allowing the infection rate to vary with time and the number of infected individuals. Such modifications enable the incorporation of external interventions, such as lockdowns or vaccination campaigns, as well as psychological and behavioral adaptations^^>\cite{buonomo2020}. 

The enriched dataset, combining real and synthetic data, is then used to train the neural networks, allowing them to effectively learn epidemic dynamics and improve predictive accuracy. By leveraging the flexibility of FNNs and the sequential modeling capabilities of NARs, our approach captures both short-term and long-term epidemic trends. Unlike PINNs, which explicitly encode governing equations into the network architecture, our method enables networks to infer dynamics implicitly through data-driven learning. This eliminates the need for embedding differential equations directly into the training process, enhancing adaptability to different scenarios while reducing computational complexity.

The remainder of this paper is structured as follows. In Section \ref{sec:prob_setting}, we recall some standard features about the classical SIR model and we introduce the social-SIR model with a saturated incidence rate used to produce synthetic data. Section \ref{sec:estimation} defines the optimization problems used to estimate parameters and the incidence function during the pre-lockdown, lockdown and post-lockdown phases in Italy and Spain. In Section \ref{sec:NN_estimation}, we present the neural network structures, distinguishing between FNNs and NARs, and we perform some detailed numerical experiments to validate our strategy. Finally, Section \ref{sec:conclusion} summarizes our findings and outlines potential future research directions.

\section{Problem setting} \label{sec:prob_setting}
The governing equations of the standard SIR model, which describes the spread of an infectious disease in a population of \(N\) individuals divided into three compartments—Susceptible (\(S\)), Infected (\(I\)), and Recovered (\(R\))—are given by  
\begin{equation}\label{eq:deterministicSIR}
	\begin{split}
		&\frac{dS(t)}{dt} = -\beta \frac{S(t) I(t)}{N},\\
		&\frac{dI(t)}{dt}  = \beta \frac{S(t) I(t)}{N} -\gamma I(t),\\
		&\frac{dR(t)}{dt} = \gamma I(t),
	\end{split}
\end{equation}  
where the total population corresponds to \(N = S(t) + I(t) + R(t)\), and the system is complemented with suitable initial conditions. The parameters \(\beta > 0\) and \(\gamma > 0\) represent the disease transmission rate and the recovery rate, respectively, and are assumed to be constant in \eqref{eq:deterministicSIR}. For simplicity, we also absorb \(N\) into the coefficient \(\beta\) in the following. 

 Despite its simplicity, the above model has been shown to be able to capture some features related to the spread of an epidemic. However, as expected, it also suffers of some drawbacks due to the inability to describe more complex structured solutions which has been observed for instance during the last COVID-19 epidemic^^>\cite{Gatto}. In particular, recent studies have highlighted the crucial role of transmission coefficients in relation to social behavior, which significantly influences pathogen spread^^>\cite{Plos}. In this research direction, the works in^^>\cite{Plos, dimarco2021kinetic} explored how controlling the number of daily social contacts can guide partial lockdown strategies. This is type of phenomena we address in this work and which cannot be captured with the system \eqref{eq:deterministicSIR}. Building on this perspective, some of the authors of the present work have developed in the recent past a new mathematical framework linking social contact distributions with disease transmission in multi-agent systems in^^>\cite{dimarco2021kinetic}. This is achieved by integrating a classical compartmental epidemiological model with kinetic equations governing social contact formation. The resulting model is given by  
\begin{equation}\label{eq:deterministicSIR_H}
	\begin{split}
		&\frac{dS(t)}{dt} = -\beta S(t) I(t) H(t,I(t)),\\
		&\frac{dI(t)}{dt}  = \beta S(t) I(t)  H(t,I(t)) -\gamma I(t),\\
		&\frac{dR(t)}{dt} = \gamma I(t),
	\end{split}
\end{equation}  
where \(H(t,I(t))\) is a macroscopic incidence rate^^>\cite{dimarco2021kinetic}, representing a time-dependent adjustment to infection rates based on behavioral and policy changes (e.g., lockdowns). In the sequel we consider the following bounded function that moderates infection spread
\begin{equation}\label{eq:H}
	H(t,I(t)) = \frac{1}{1+\frac{\nu}{\beta} \int_0^t I(s) ds}, 
\end{equation}  
where \(\nu > 0\) is a constant that it has been shown to reflect the variance in contact rates within the population. The system \eqref{eq:deterministicSIR_H} even if built on a different basis, shares conceptual similarities with the models introduced in^^>\cite{capasso1978}, where a nonlinear, bounded function \(g(I)\) modulates disease transmission. This approach accounts for the fact that, as the number of infected individuals increases, the effective infection force may decrease due to factors such as behavioral adaptations or resource limitations. On the same direction in^^>\cite{liu1986influence} it has been introduced a nonlinear, saturated incidence rate of the form  
\begin{equation}\label{sat_incodencerate}
	g(I)S = \frac{ \beta I^{m}S}{1 + \alpha I^{k}}, \quad \text{with} \quad \frac{dS(t)}{dt} = -g(I(t)) S(t)
\end{equation}  
where \(\beta I^{m}S\) represents the baseline infection force, and \(\frac{1}{1 + \alpha I^{k}}\) accounts for behavioral responses, such as increased caution among susceptible individuals. Here, \(\beta, m, k\) are positive constants, and \(\alpha\) is a non-negative constant^^>\cite{cui2017complex, tang2008coexistence, ruan2003dynamical}.  A specific case, with \(m = 1, k = 1\), yields the form proposed in^^>\cite{capasso1978} to model the crowding effect observed during the 1973 cholera outbreak in Bari: 
\begin{equation}\label{sat_incodencerate1}
	g(I) = \frac{ \beta I}{1 + \alpha I}.
\end{equation}  
Alternatively, setting \(m = 1, k = 2\) results in the Monod–Haldane (M–H) function,  
\begin{equation}\label{sat_incodencerate2}
	g(I) = \frac{ \beta I}{1 + \alpha I^{2}},
\end{equation}  
which is commonly used to model physiological constraints in infection dynamics.  

The so-called social-SIR model \eqref{eq:deterministicSIR_H} exhibits a similar incidence rate behavior, extending the form in \eqref{sat_incodencerate1} to account for the impact of protective measures such as lockdown strategies and infection control. Unlike classical models, this generalized framework incorporates a time-dependent, saturated incidence rate to capture the cumulative effects of infection dynamics and behavioral adaptations in response to an epidemic. This model has previously been successfully employed to analyze the effects of partial lockdown strategies in various European countries. For the purpose of demonstrating the proposed data augmentation strategy, we select a dataset corresponding to a period marked by COVID-19 lockdowns. It is worth to mention that, the complexity introduced by lockdown measures, makes the forecasting process more challenging compared to standard scenarios. Therefore, in the following, we employ system \eqref{eq:deterministicSIR_H} to develop our data-driven strategy and the associated deep learning framework, enabling epidemic forecasting.

\section{On the the construction of a data driven social-SIR model}\label{sec:estimation}
In the context of epidemic prediction, mathematical models such as those presented in the previous section provide valuable insights into the spread of infectious diseases. However, these models often rely on predefined assumptions and parameters that may not fully capture the complexities of real-world epidemics. Neural networks offer a powerful alternative by learning directly from data, enabling them to model complex, nonlinear relationships that traditional epidemic models may struggle to represent.  
However, one major limitation of training neural networks exclusively on available data is that, while they can effectively interpolate within a given time interval, they struggle to make accurate predictions beyond the observed data points^^>\cite{la2020epidemiological,amendolara2023lstm,bertaglia2022asymptotic}. To address this issue, we propose an approach that combines a mathematical epidemiological model—the social-SIR model \eqref{eq:deterministicSIR_H}—with real outbreak data. The key idea is to generate additional synthetic information by solving this generalized SIR model, thereby filling the gaps between available real or measured data points.  
With this enriched dataset, one can train a neural network more effectively to predict future infection trends in a given population. As we demonstrate later, this approach leverages both real and synthetic data, significantly improving predictive accuracy and enhancing the model’s ability to forecast epidemic dynamics. However, to implement this strategy, the model must be capable of interpolating through the experimental data points. Therefore, in the following section, we discuss how the parameters in \eqref{eq:deterministicSIR_H} are defined.

\subsection{Parameters estimation}\label{sec:params_estimation}
To develop our approach, we focus on the COVID-19 outbreak in two European countries, Italy and Spain, during the 2020 pandemic.  Recognizing that data collection processes inherently introduce errors, particularly due to underreporting of number of infected individuals in the pandemic’s early stages, we assume an intrinsic uncertainty is embedded in the collected data, capturing the variability and limitations of real-world measurements, \cite{chowell2017fitting}. While the development of advanced uncertainty quantification techniques falls outside the scope of this work, it represents a promising avenue for future research.  In the present study, we instead adopt the well?established framework described in \cite{dimarco2021kinetic, albi2022kinetic}.
Given that both countries implemented restrictive measures during this period, we proceed as follows. We define two optimization phases: in the first phase, we estimate the parameters \(\beta\) and \(\gamma\) for the pre-lockdown period using the model \eqref{eq:deterministicSIR_H} with $H(t,I(t))=1$; in the second phase, we estimate the function \(H(t,I(t))\) from lockdown and post-lockdown data. 

The COVID-19 data for Italy are available from the
\href{https://github.com/pcm-dpc/COVID-19/tree/master/dati-andamento-nazionale}{GitHub repository } of the Italian Civil Protection Department, while data for Spain can be accessed \href{https://github.com/CSSEGISandData/COVID-19}{here}. 
In practice, to estimate the parameters $\beta$ and $\gamma$, we formulate the following minimization problem:  
\begin{equation}\label{eq:min_beta_gamma}
	\min_{\beta,\gamma} \mathcal{J}(\hat{I},I),
\end{equation}  
subject to the constraint that the equations of system \eqref{eq:deterministicSIR_H} have to be satisfied when $H(t,I(t))=1$. The objective function to be minimized is defined as  
\begin{equation}\label{eq:functional}
	\mathcal{J}(\hat{I}, I) =  \|\hat{I}(t) - I(t)\|^2,
\end{equation}  
where \(\hat{I}(t)\) represents the reported infected cases from the data set at disposal, while \(I(t)\) is the model output, and \(\Vert \cdot \Vert\) denotes the \(L_2\) norm.  Throughout the manuscript, the variables $S$, $I$, and $R$ and the data are assumed to be normalized and to denote fractions of the total population.

\begin{figure}[H] 
	\centering
	\includegraphics[width=0.458\linewidth]{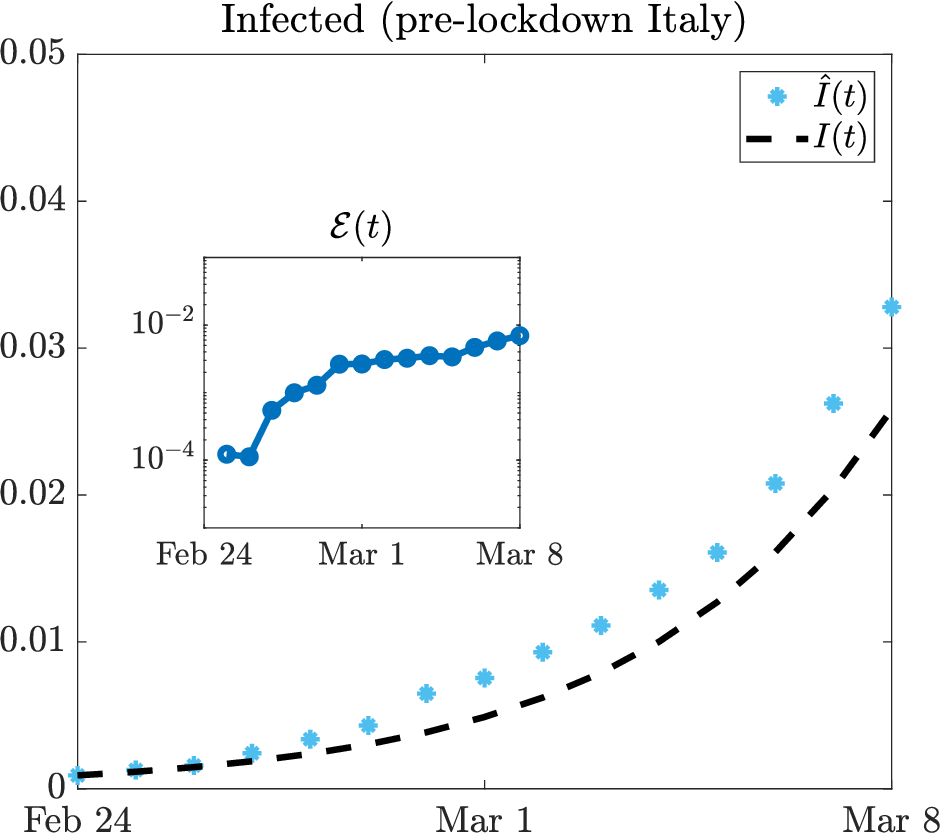}
	\includegraphics[width=0.458\linewidth]{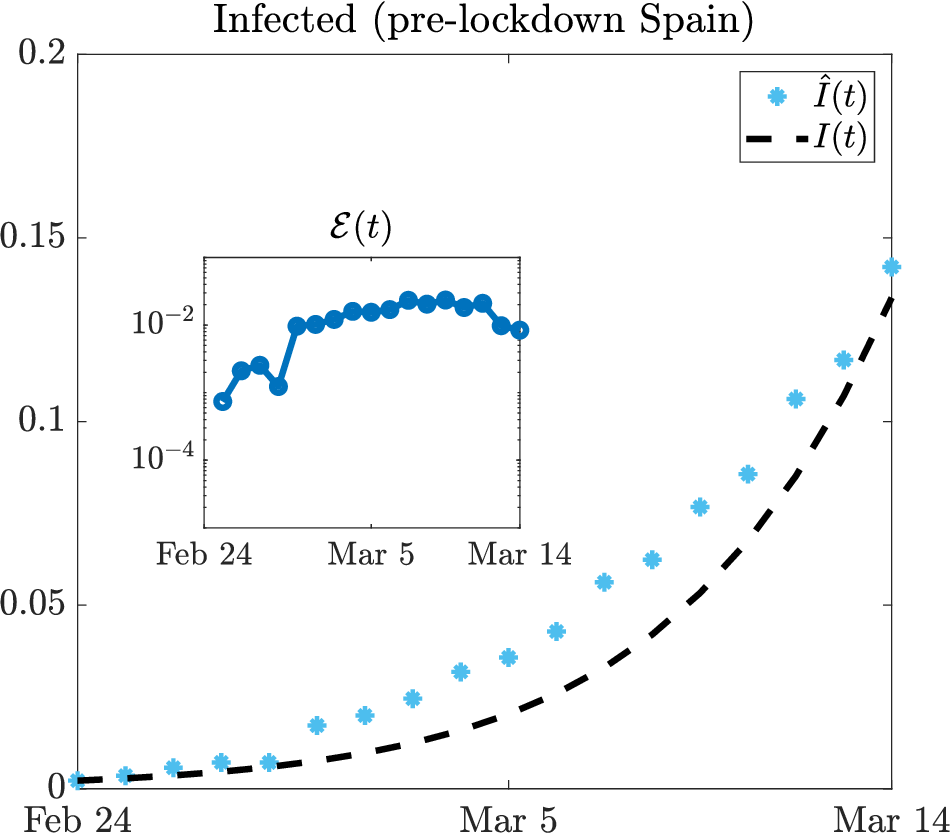}
	\caption{Pre-lockdown phase. Dynamics of the infected population computed  by solving the social-SIR model \eqref{eq:deterministicSIR_H} with $\beta$ and $\gamma$ as in \eqref{eq:beta_gamma_est} (Italy), and as in  \eqref{eq:beta_gamma_est_spain} (Spain) and $H(t,I(t))=1$, compared with the experimental data. In the embedded plot the pointwise error between the two quantities computed as in \eqref{eq:error}. On the left, Italy. On the right, Spain.}
	\label{fig:pre_lockdown}
\end{figure}

We consider data over the time interval \([t_0,t_l)\), where \(t_0\) (set to February 24, 2020) marks the initial time, and \(t_l\) represents the lockdown start date (March 9 for Italy and March 14 for Spain). We search for a minimum under the constraints  
\begin{equation}\label{eq:constraints_beta_gamma}
	0\leq \beta \leq 0.6,\qquad 0\leq \gamma \leq 0.06,
\end{equation}  
where the bounds are chosen accordingly to the available estimates^^>\cite{Gatto}. To solve the optimization problems in \eqref{eq:min_beta_gamma}, we rely on the use of the Matlab functions \texttt{fmincon} in combination with a RK4 integration method of the system of ODEs.  The result of the optimization procedure on the pre-lockdown phase gives, for Italy, the following values 
\begin{equation}\label{eq:beta_gamma_est}
	\beta_{est}^{It} = 0.3107, \qquad \gamma_{est}^{It} = 0.0328,
\end{equation}
while for Spain 
\begin{equation}\label{eq:beta_gamma_est_spain}
	\beta_{est}^{Sp} =  0.2340, \qquad \gamma_{est}^{Sp} = 0.0224.
\end{equation}
Both the estimated transmission rates $\beta_{est}$ and the recovery rates $\gamma_{est}$ have units of inverse time (e.g., $	\texttt{day}^{-1}$), ensuring the correct dimensional consistency of the model equations.
In Figure \ref{fig:pre_lockdown} on the left the evolution of the dynamics in the pre-lockdown phase in Italy obtained by solving the social-SIR model \eqref{eq:deterministicSIR_H} with $H(t,I(t))=1$, using the estimated parameters $\beta = \beta_{est}^{It}$ and $\gamma = \gamma_{est}^{It}$ in \eqref{eq:beta_gamma_est}. This is compared, in the same image, to the experimental data from the pre-lockdown phase, spanning up to time $t_l=14$ with a time step $\Delta t=1$ corresponding to $14$ days with $1$ day interval between the measurements. The estimated dynamics aligns quite well with the available data, as confirmed by the embedded plot which shows the pointwise error over time between the solution of the generalized SIR model and the data at disposal, i.e. 
\begin{equation}\label{eq:error}
	\mathcal{E}(t) = \vert \hat{I}(t) - I(t)\vert.
\end{equation}
In the same Figure on the right the pre-lockdown dynamics in Spain computed by solving the social-SIR model \eqref{eq:deterministicSIR_H}  with $H(t,I(t))=1$, using the estimated parameters $\beta = \beta_{est}^{Sp}$ and $\gamma = \gamma_{est}^{Sp}$ in \eqref{eq:beta_gamma_est_spain} compared with the experimental data spanning up to time $t_l = 18$. The two solutions are in good agreement as confirmed by the plot showing the pointwise error computed as in \eqref{eq:error}. 

\begin{figure}[H]
	\centering
	\includegraphics[width=0.458\linewidth]{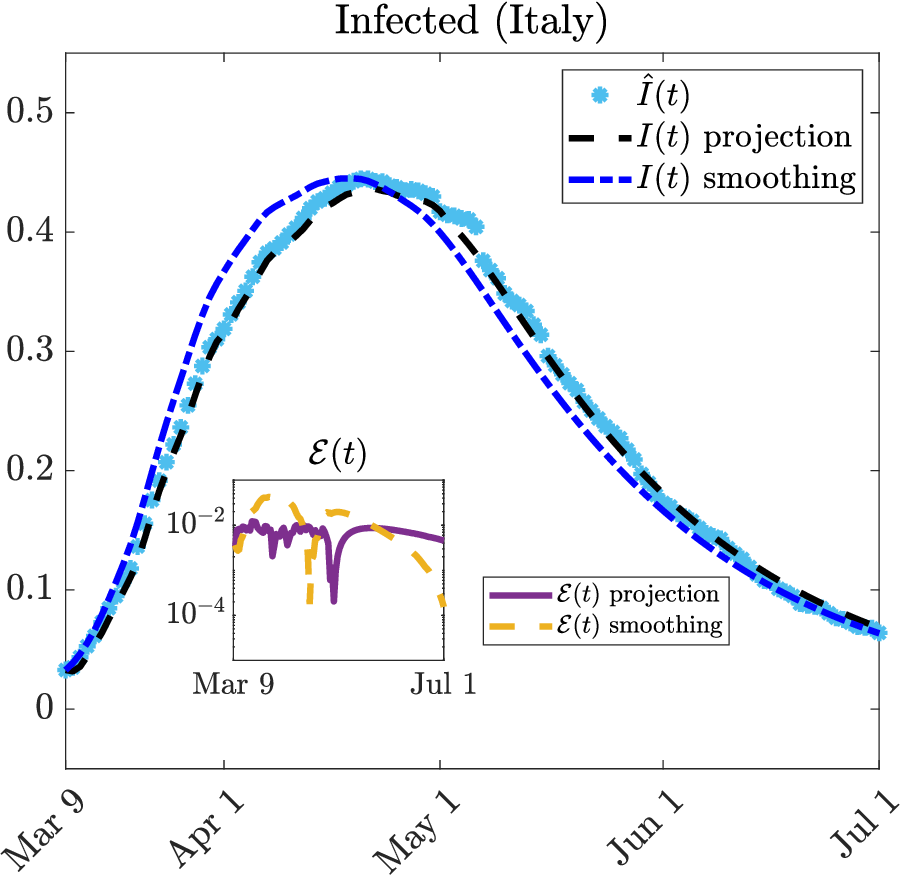}
	\includegraphics[width=0.458\linewidth]{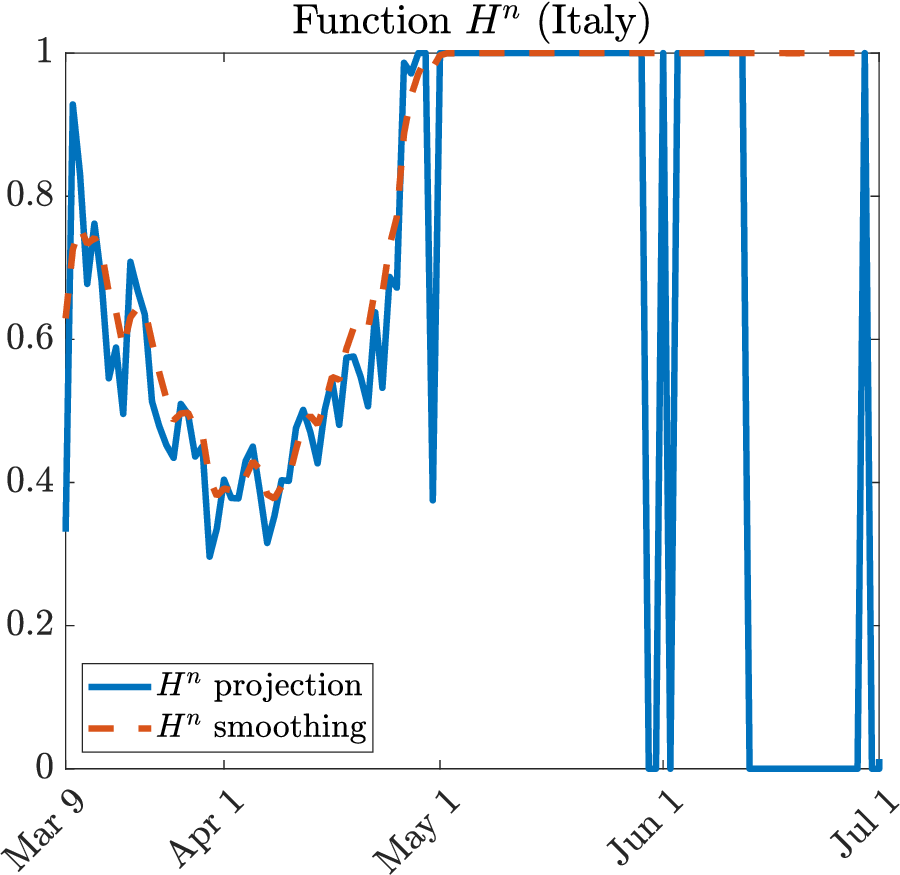}\\
	\includegraphics[width=0.458\linewidth]{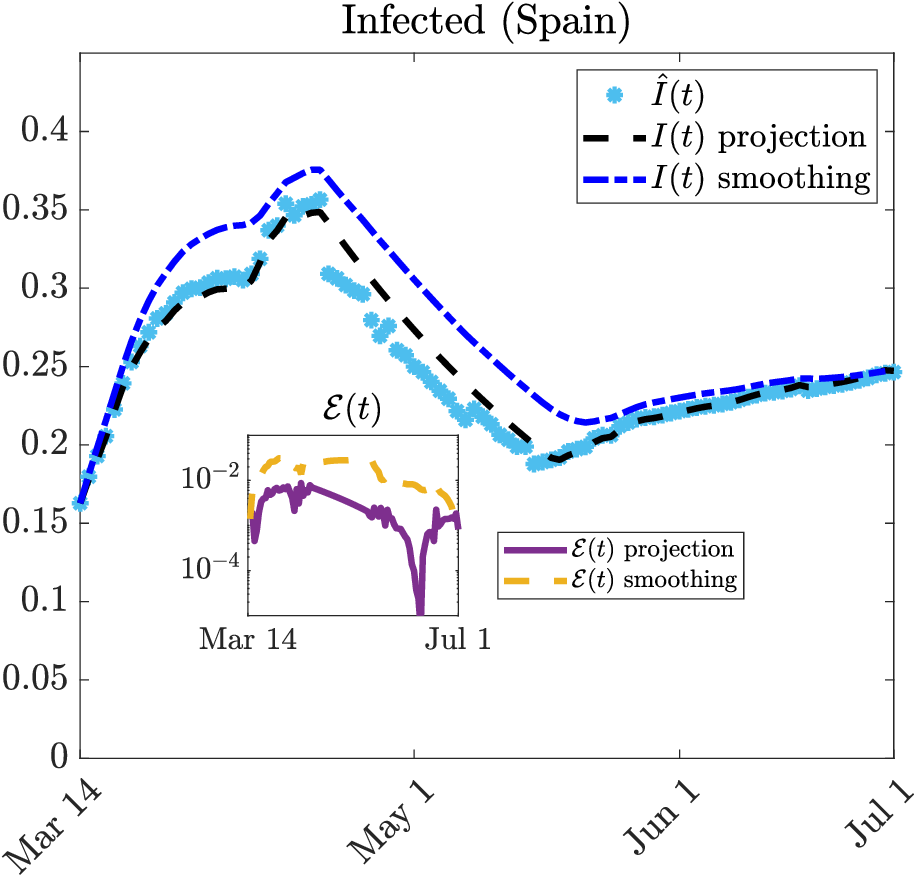}
	\includegraphics[width=0.458\linewidth]{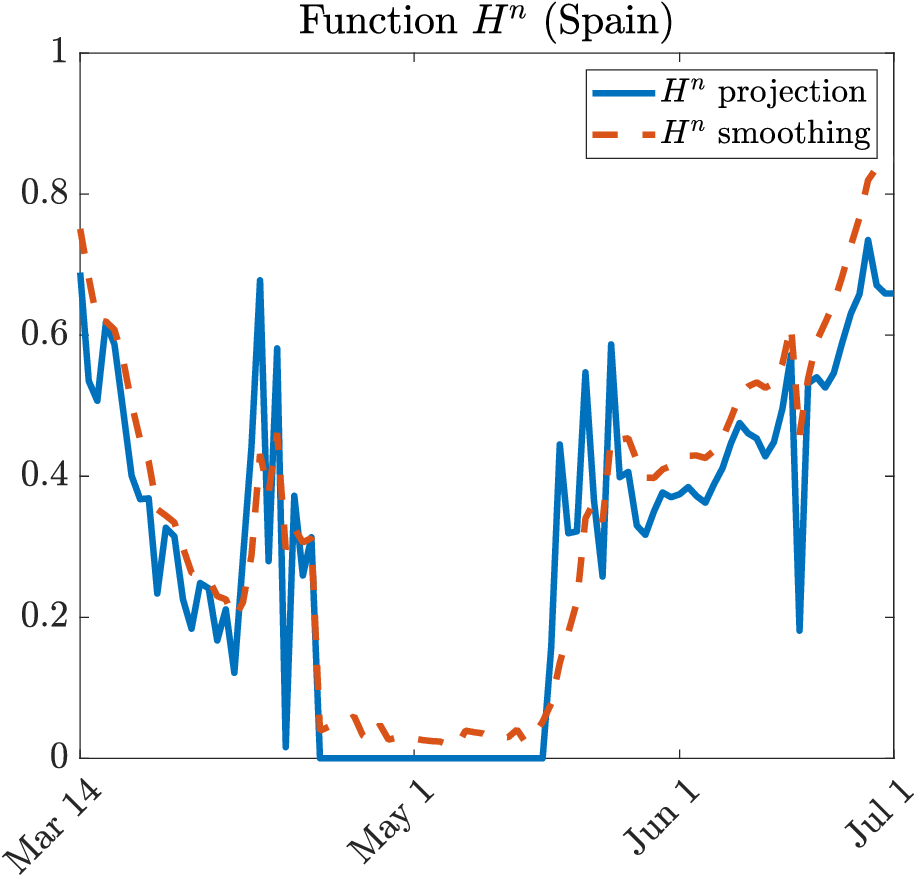}
	\caption{Lockdown and post-lockdown phase. On the left the dynamics of the infected population computed by solving the generalized SIR model \eqref{eq:deterministicSIR_H} with $\beta$ and $\gamma$ as in \eqref{eq:beta_gamma_est} (Italy), and as in \eqref{eq:beta_gamma_est_spain} (Spain), assuming the contact function $H$ to be defined as in \eqref{eq:H_instaneous} and, alternatively, with the smoothed version \eqref{eq:smoothing_projection} setting $k=1.6$ compared with the experimental data. In the embedded plot the pointwise error between the two infected populations computed as in \eqref{eq:error}.  On the right the contact function $H$ computed as in \eqref{eq:H_instaneous} and its smoothed version computed introducing an operator as in \eqref{eq:smoothing_projection}. First row, Italy. Second row, Spain.}
	\label{fig:post_lockdown}
\end{figure}

Once \(\beta\) and \(\gamma\) have been estimated, we proceed to determine the pointwise values \(H(t,I(t))\) in \eqref{eq:deterministicSIR_H} in such a way that the model fits as well as possible the lockdown and post-lockdown data, i.e., data in the interval \([t_l,T]\), where \(T\) (set to July 1st) is the final observation time for both countries. To this end, we solve the following new optimization problem  
\begin{equation}\label{eq:min_H}
	\min_{H(t)} \mathcal{J}(\hat{I},I),
\end{equation}  
subject to the constraint that the solution has to satisfy the equations of the system \eqref{eq:deterministicSIR_H}, with \(\mathcal{J}(\cdot)\) defined as in \eqref{eq:functional}, and that 
\begin{equation}\label{eq:constraint_H}
	0\leq H(t) \leq 1.  
\end{equation}
In order to determine the optimal value of $H$, we split the time interval $[t_l, T]$ into $N_t$ time intervals of length $\Delta t$, and we solve a sequence of optimal control problems
\begin{equation}\label{eq:min_H_discr}
	\min_{H} \mathcal{J}^{\Delta t}(\hat{I},I) := \int_t^{t+\Delta t} \Vert \hat{I}(s)-I(s)\Vert^2 ds,
\end{equation}
s.t. the explicit time discretization of the generalized SIR dynamics
\begin{equation}\label{eq:discr_eqs} 
	\begin{split}
		&S^{n+1} = S^{n} + \Delta t \beta I^n S^n H^n,\\ 
		&I^{n+1} = I^{n} + \Delta t \left(\beta I^n S^n H^n -\gamma I^n\right),\\
		&R^{n+1} = R^n + \Delta t  \gamma I^n,
	\end{split}
\end{equation}
where we denoted by $t^n = n \Delta t$ for any $n=1,\ldots,N_t$. We then consider a discretization of the integral in \eqref{eq:min_H_discr} using the rectangle rule  
\begin{equation}\label{eq:discr_functional}
	\mathcal{J}^{\Delta t}(\hat{I},I) = \Delta t \left( \hat{I}^{n+1}-I^{n+1}\right).
\end{equation}
By direct computation, it can be shown that the expression for $H$ at each discrete time $t^n$ is given by  
\begin{equation}\label{eq:H_instaneous}
	H^n = \mathbb{P}_{[0,1]}\left( \frac{\hat{I}^n - I^n + I^n \gamma \Delta t}{I^n S^n \beta \Delta t}\right),
\end{equation}
where $\mathbb{P}_{[0,1]}(\cdot)$ denotes the projection over the interval $[0,1]$. 
As demonstrate later, the function $H$ computed as in \eqref{eq:H_instaneous}  tends to exhibit oscillatory behavior over time. To reduce these oscillations, we propose to replace the projection operator $\mathbb{P}_{[0,1]}(\cdot)$ in equation \eqref{eq:H_instaneous} with a smoother approximation of the type 
	\begin{equation}\label{eq:smoothing_projection}
		\mathbb{S}(x) = \frac{1}{1+e^{-kx}},
	\end{equation}
	for a certain $k>0$, and $x\in \mathbb{R}$. 
	This modification helps to regularize the value of $H$.
The results of this study are reported in Figure \ref{fig:post_lockdown} where on the left we show the lockdown and post-lockdown experimental data and the solution of the generalized SIR model with the contact function $H$ as in \eqref{eq:H_instaneous} and, alternatively, with the smoothed version \eqref{eq:smoothing_projection} setting $k=1.6$. The estimated and true solutions are in good agreement as shown in the embedded plot which represents the corresponding pointwise error computed as in \eqref{eq:error}. On the right, we plot the estimated values of $H^n$ over time, comparing the values obtained from the original projection-based approach and the smoothed operator. The function $H$ computed as in \eqref{eq:H_instaneous} exhibits oscillations over time, suggesting potential measurement errors in the data such as a non-exponential decay of the number of infected especially in the late post-lockdown phase, mimicking the intrinsic uncertainty that may affect data. These oscillations are significantly reduced when using the smoothed operator.  However, a reduction in accuracy of the numerical solution can be observed. First row, Italy. Second row, Spain.

\section{Estimation and prediction using real COVID-19 data and Neural Networks design for epidemic prediction}\label{sec:NN_estimation} 
In this section, we consider two distinct types of neural network architectures: Feed-Forward Neural Networks (FFNs) and Nonlinear Autoregressive Networks (NARs). Both approaches are implemented using a general neural network framework with $L+1$ layers, as follows
\begin{equation}
	\begin{split}
		&x^1 = W^1 x + b^1, \qquad x^{l} = \sigma \circ (W^lx^{l-1} + b^l), \quad l=2,\ldots,L-1,\\
		&U^{NN}(x;{\bf b},{\bf W}) = W^L z^{L-1}+ b^L,
	\end{split}
\end{equation} 
where $x = (x^1,\ldots,x^d) \in \mathbb{R}^d$ is the input, $\sigma$ is an activation function  to convert the input signals to output signals, and ${\bf W} = (W^1,\ldots, W^L)$ and ${\bf b} = (b^1,\ldots,b^L)$ represent the weights and biases. The weights and biases associated with these connections serve as the parameters of the network, adjusted iteratively through techniques like gradient descent or Adam method during the training phase to minimize prediction errors and enhance model performance. In our context, the idea is to solve the following minimization problem
\begin{equation}\label{eq:min_prob}
	\theta^* = \arg \min_{\theta} \mathcal{L}(\theta),
\end{equation}
to find the optimal parameters $ \theta = (W^1,b^1,\ldots,W^L,b^L)$ which minimize a certain loss function $\mathcal{L}(\theta)$  that includes errors from both real and synthetic data generated by solving the generalized SIR model. Real data provides an accurate representation of the epidemic dynamics, while synthetic data fills gaps between observed data points and provides high-resolution trajectories. 

The loss function is expressed as
\begin{equation}\label{eq:combined_loss}
	\mathcal{L}(\theta) = \omega_r \mathcal{L}_r(\theta) + \omega_s \mathcal{L}_s(\theta),
\end{equation}
where $\omega_r, \omega_s>0$ are given weights, and 
\begin{equation}
	\mathcal{L}_j(\theta) =\frac{1}{N_j} \sum_{i=1}^{N_j} \|  I_j (t) - I_j^{NN}(t, \theta) \|^2,
\end{equation}
for $j=\{r,s\}$,  where $I_j(t)$ is the actual number of infected. These are either the real data or the data generated by solving the social SIR model (synthetic data), while $I_j^{NN}(t, \theta)$  is the predicted state. Finally, $N_r$, $N_s$ refer to the number of available real and synthetic data points respectively. With this choice, it is possible to appropriately weight the contribution of real data alongside the data generated by the model.
\paragraph{Feed-Forward Neural Networks.}
FFNs are a foundational type of neural network architecture designed to map inputs directly to outputs without considering temporal dependencies. For epidemic prediction, we design a FFN that takes the current time 
$t$ as input and produces the value of infected  
$I(t)$ as output. This approach relies on learning a direct mapping between the time and the epidemic dynamics, effectively interpolating the trajectory based on the training data.
The architecture of the FFN consists of the following components
\begin{itemize}
	\item \textbf{Input Layer:} The single input feature is the current time $t$.
	\item \textbf{Hidden Layers:} One or more fully connected layers with nonlinear activation functions (e.g., ReLU or tanh) capture the relationship between time and the epidemic states.
	\item \textbf{Output Layer:} The output consists of a node corresponding to the predicted value of the infected at time $t$, i.e.  $I(t)$. 
\end{itemize}
It is well known that, while FFNs are effective for interpolation tasks, their inability to model sequential dependencies limits their accuracy in predicting future epidemic dynamics. This limitation arises because FFNs lack memory mechanisms, meaning they cannot leverage information from previous states to influence current predictions. For this reason, in the following we propose an alternative strategy.
\paragraph{Nonlinear Autoregressive Networks.}
To capture the temporal dynamics of the epidemic, we then rely on Nonlinear Autoregressive Networks (NARs) capable of modeling sequential data. Unlike FFNs, NARs incorporate information about past states, making them particularly effective for time-series predictions.  The NAR is designed to predict the future state of the epidemic based on the previous state of the solution. The architecture consists of the following components:
\begin{itemize}
	\item \textbf{Input Layer:} The input features are the values of the infected at time $t-d,\ldots,t-1$, i.e. $I(t-d),\ldots, I(t-1)$ for a certain $d\geq1$ .
	\item \textbf{Hidden Layers:} One or more fully connected layers with nonlinear activation functions (e.g., ReLU or tanh) capture nonlinear relationships in the time series data. 
	\item \textbf{Output Layer:} Produces the predicted next state, $I(t)$, representing the system's state at time $t$.
\end{itemize}
Once trained, NAR networks can predict the future states of the epidemic compartments by recursively feeding its outputs back as inputs for subsequent time steps (closed loop strategy). This enables the generation of full epidemic trajectories from a given initial condition.

\subsection{Neural network for epidemic prediction} \label{sec:numerical_experiments}
\begin{figure}[H]
	\centering
	\includegraphics[width=0.458\linewidth]{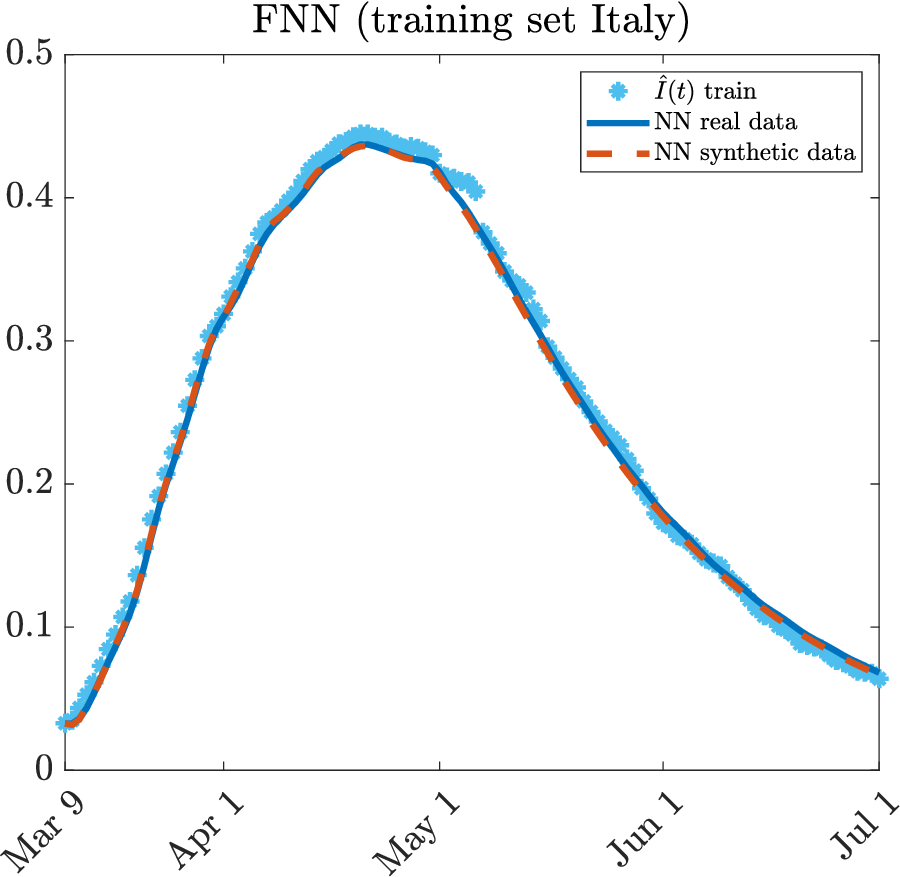}
	\includegraphics[width=0.459\linewidth]{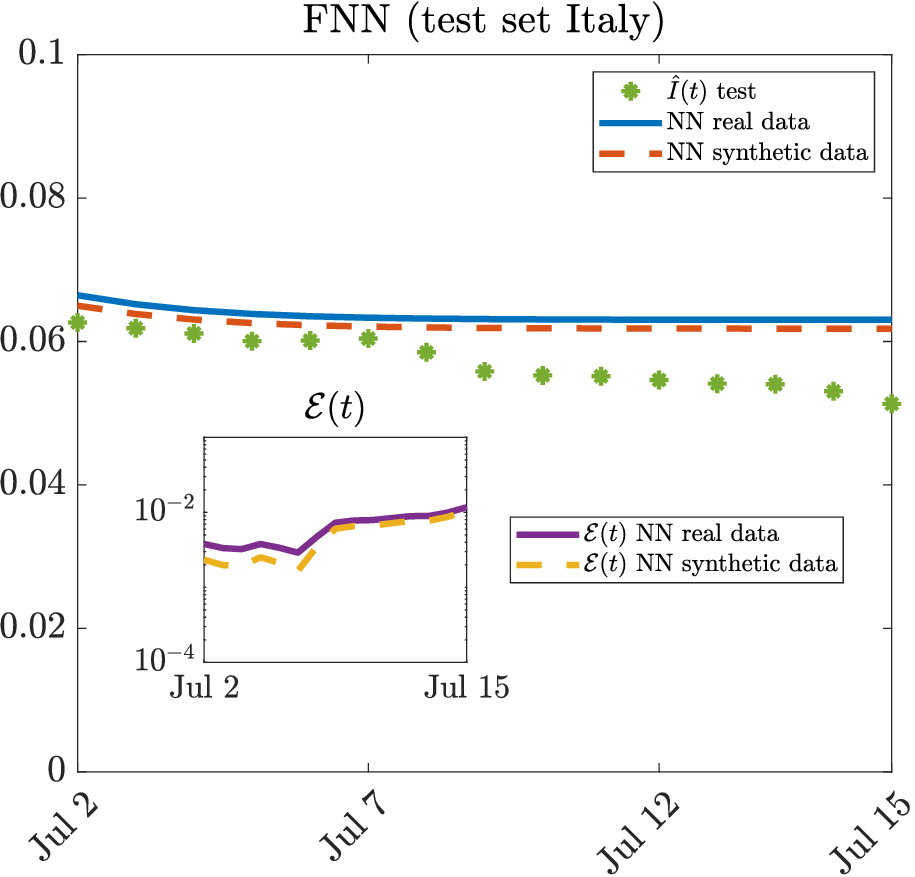}
	\caption{Feedforward neural network for Italy. On the left the dynamics of the infected population computed by training a FNN on real and synthetic data and tested over the train interval $[t_l,T]$ compared to the available data. On the right the solution computed with the same FNN tested over the time interval $(T,T_{test}]$ compared to the available data. Embedded the error between the available data test and the neural network solution computed as in \eqref{eq:error}.}
	\label{fig:FNN_I1}
\end{figure}

\begin{figure}[H]
	\centering
	\includegraphics[width=0.458\linewidth]{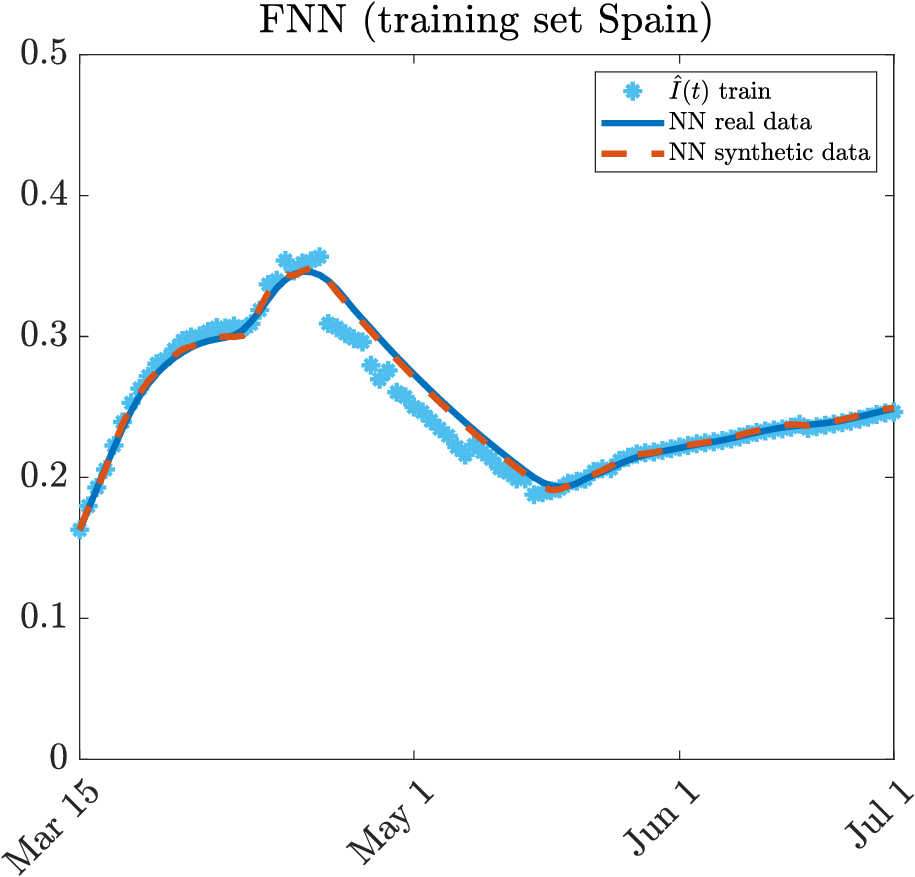}
	\includegraphics[width=0.458\linewidth]{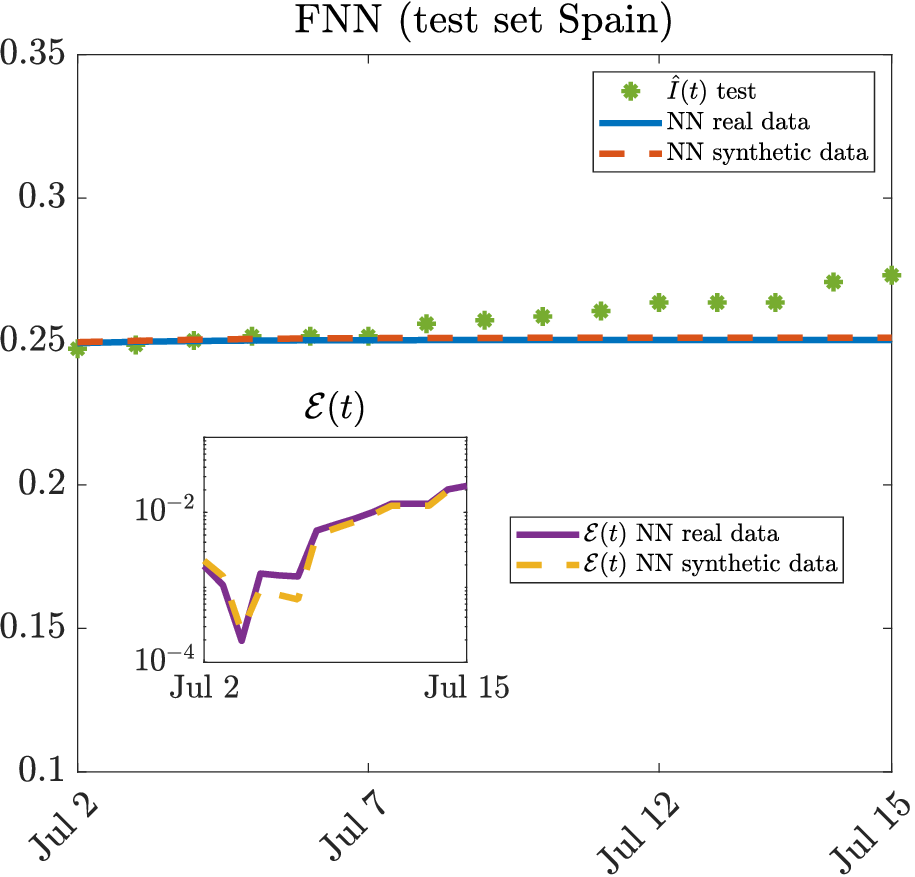}
	\caption{Feedforward neural network for the Spain case. On the left the dynamics of the infected population computed by training a FNN on real and synthetic data and tested over the train interval $[t_l,T]$ compared to the available data. On the right the solution computed with the same FNN tested over the time interval $(T,T_{test}]$ compared to the available data. Embedded the error between the available data test and the neural network solution computed as in \eqref{eq:error}. First row, Italy. Second row, Spain.}
	\label{fig:FNN_I2}
\end{figure}

To evaluate the performances of FFNs and NARs in epidemic prediction, we conduct some numerical experiments using both real and synthetic data. The real data covers the period from the lockdown start date (in Italy March 9th, in Spain March 14th) to July 1st. These data are available on a daily basis, corresponding to a time step of $\Delta t_r = 1$. The synthetic data are generated by solving the generalized SIR model given in \eqref{eq:deterministicSIR_H}, with $\beta  $ and $\gamma $ as in \eqref{eq:beta_gamma_est} (Italy), and as in \eqref{eq:beta_gamma_est_spain} (Spain) over the same time interval, using a finer time step $\Delta t_s = 0.2$, and with $H$ estimated by solving  the minimization problem in \eqref{eq:min_H}. In the following, we will focus on the solution generated by solving the social-SIR model with $H$ as in \eqref{eq:H_instaneous} to generate synthetic data. The version with the smoothing operator is not considered further due to the loss of accuracy it introduces. Furthermore, since $H$ is originally defined at daily intervals as in \eqref{eq:H_instaneous}, we use linear interpolation to extend its definition to the finer grid. 
Both networks are trained using the Levenberg-Marquardt optimizer for $500$ epochs. In the case of real data we set $\omega_r=1$ and $\omega_s=0$, while in the case of synthetic data we set $\omega_r=0$ and $\omega_s=1$ in \eqref{eq:combined_loss}.
\paragraph{Feedforward neural network.}
We design a FFN with three fully connected layers, each consisting of $10$ nodes, and use the \texttt{tanh} activation function, relying on the use of the Matlab function \texttt{feedforwardnet}.  The input to the network is the time interval $[t_l,T]$, being $t_l$ the lockdown starting date and $T$ the final time supposed to be July 1st. 
The output is the normalized number of infected $I(t)$, which corresponds to either real data or synthetic data generated as described above. We train the neural network using the Matlab function \texttt{train}. To evaluate the neural network's performance, we first compute its solution over the training set for both Italy and Spain. Specifically, we consider the time interval $[t_l,T]$ with a time step $\Delta t_{r}=1$. In the left panel of Figure \ref{fig:FNN_I1} and of Figure \ref{fig:FNN_I2}, we compare the solution obtained using a neural network trained on both real and synthetic data, with the available data for respectively the Italian and the Spanish case. The neural network provides a good approximation of the solution, making it well-suited for interpolation tasks. Next, to assess its predictive capabilities, we focus on the time interval 
$(T,T_{test}]$, where $T$ (July 1st) is the starting time and $T_{test}$ (July 15th) is the final time, using the same daily discretization. The right panels of Figures \ref{fig:FNN_I1} and \ref{fig:FNN_I2} present a comparison between the network's predictions for the test set and the daily real data available in that period. Additionally, the figures display the pointwise error between the neural network's outputs and the real data, calculated as in \eqref{eq:error}. The standard FNN as expected faces challenges in forecasting future dynamics and the data augmentation does not provide a clear improvement in the results.

\begin{figure}[htb]
	\centering
	\includegraphics[width=0.458\linewidth]{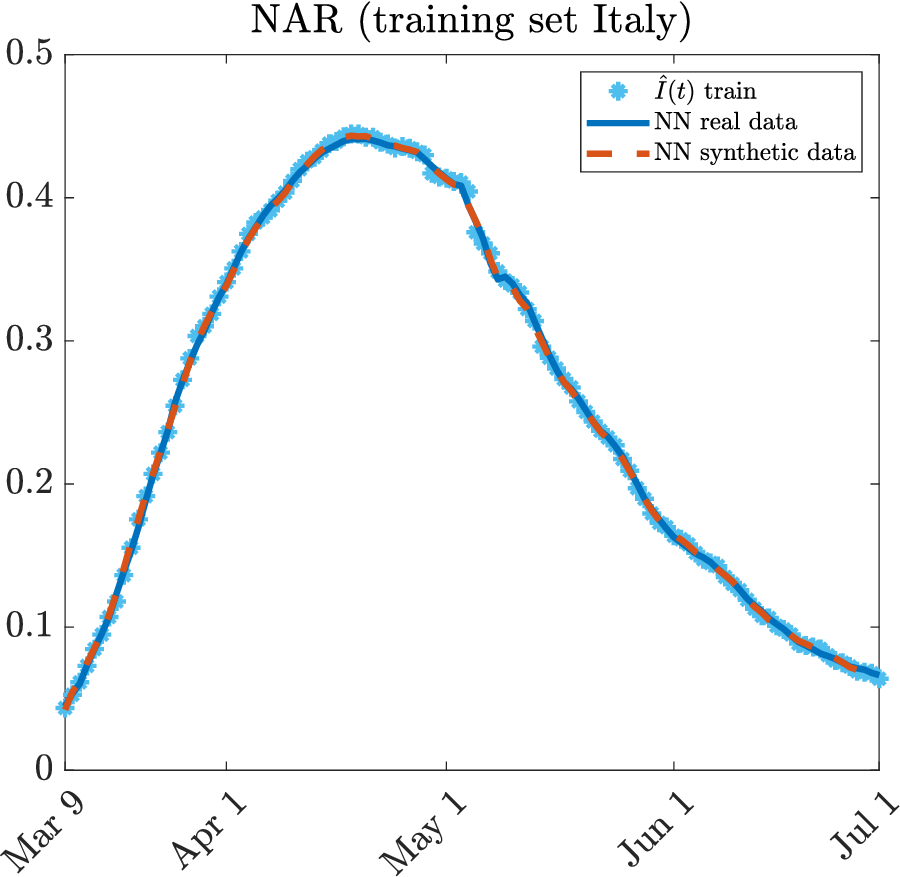}
	\includegraphics[width=0.458\linewidth]{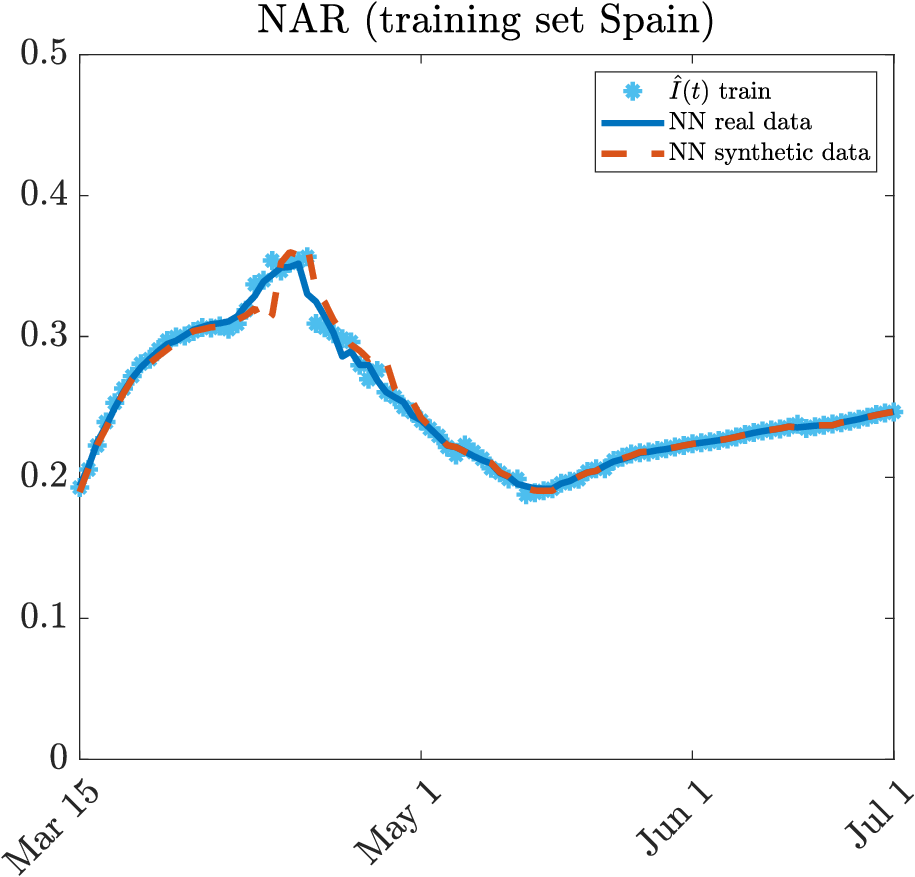} 
	\caption{Nonlinear autoregressive network: training data. Solution obtained by training a NAR network on both real and synthetic data, then testing it on the same dataset and comparing it to the available data. On the left, the case for Italy, and on the right, the one for Spain.}
	\label{fig:NAR_training}
\end{figure}

\begin{figure}[H]
	\centering
	\includegraphics[width=0.458\linewidth]{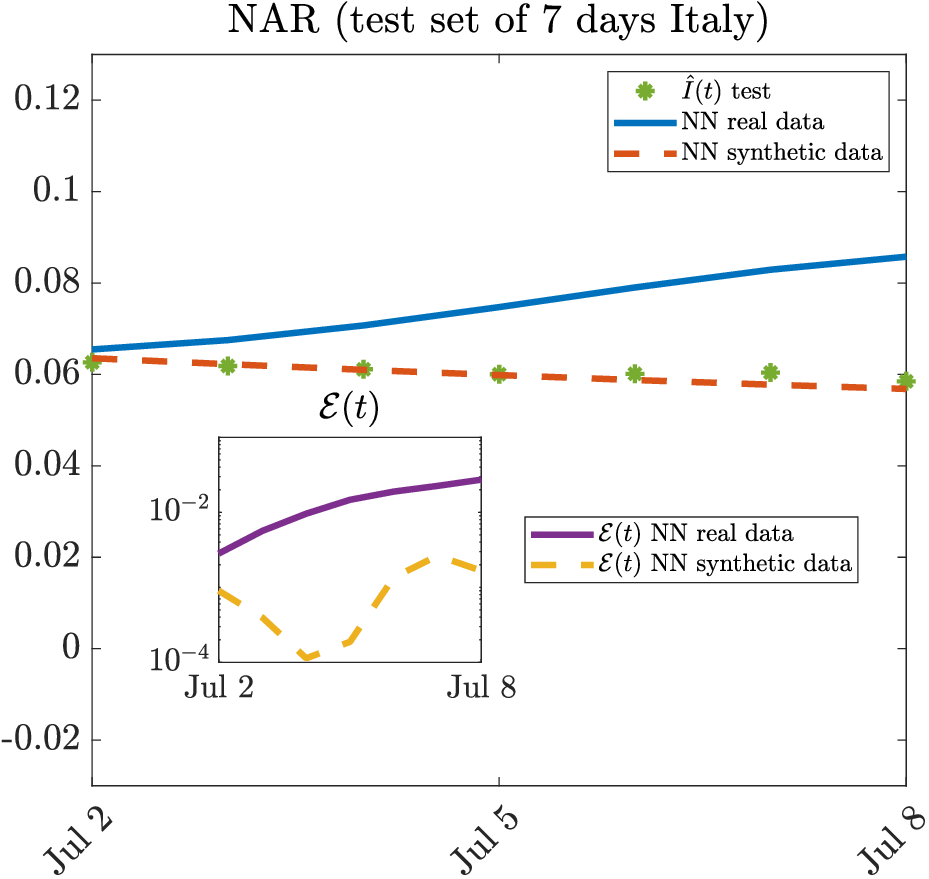}
	\includegraphics[width=0.458\linewidth]{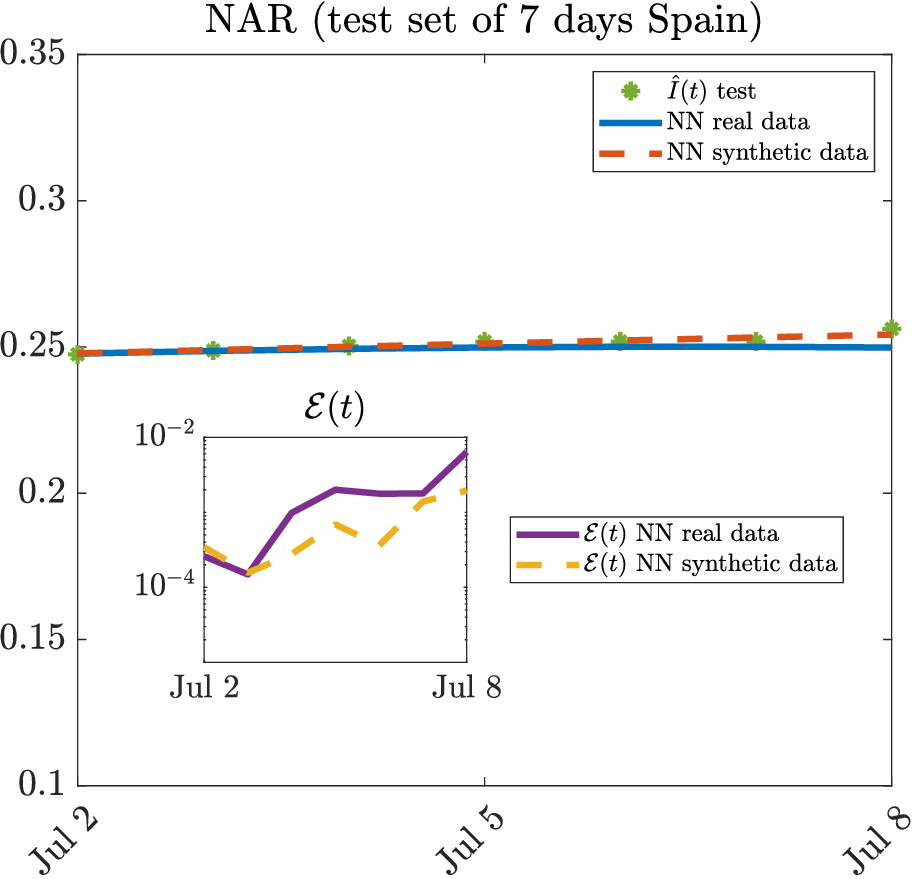}
	\includegraphics[width=0.458\linewidth]{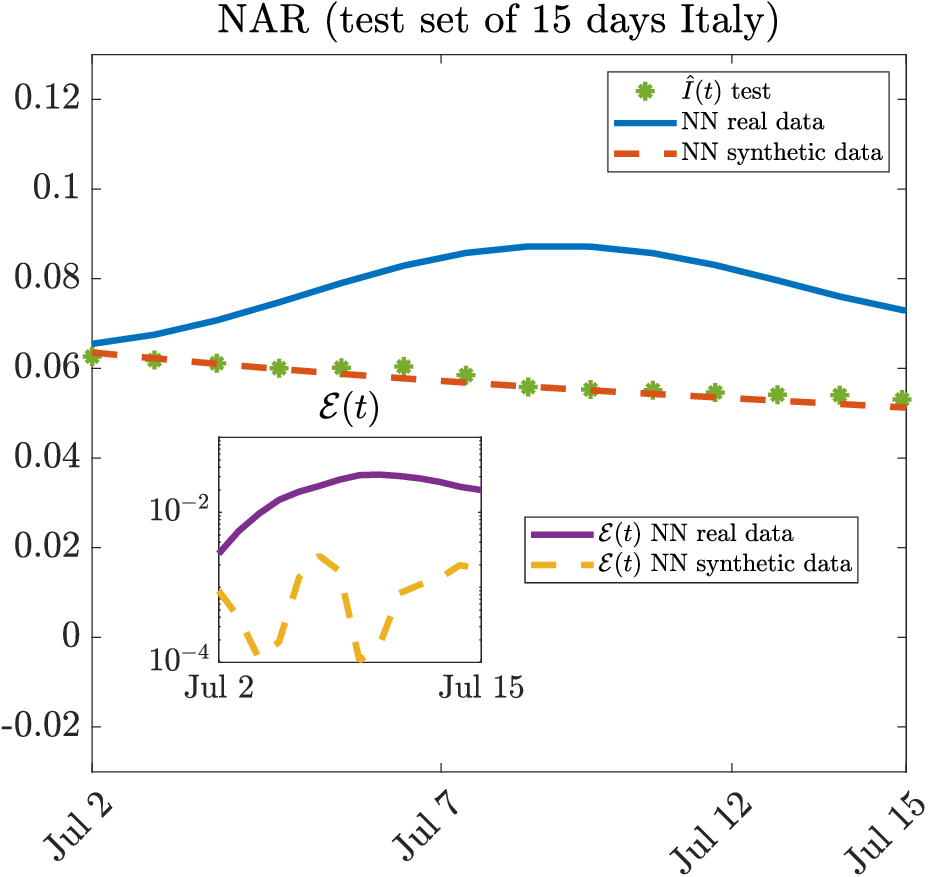}
	\includegraphics[width=0.458\linewidth]{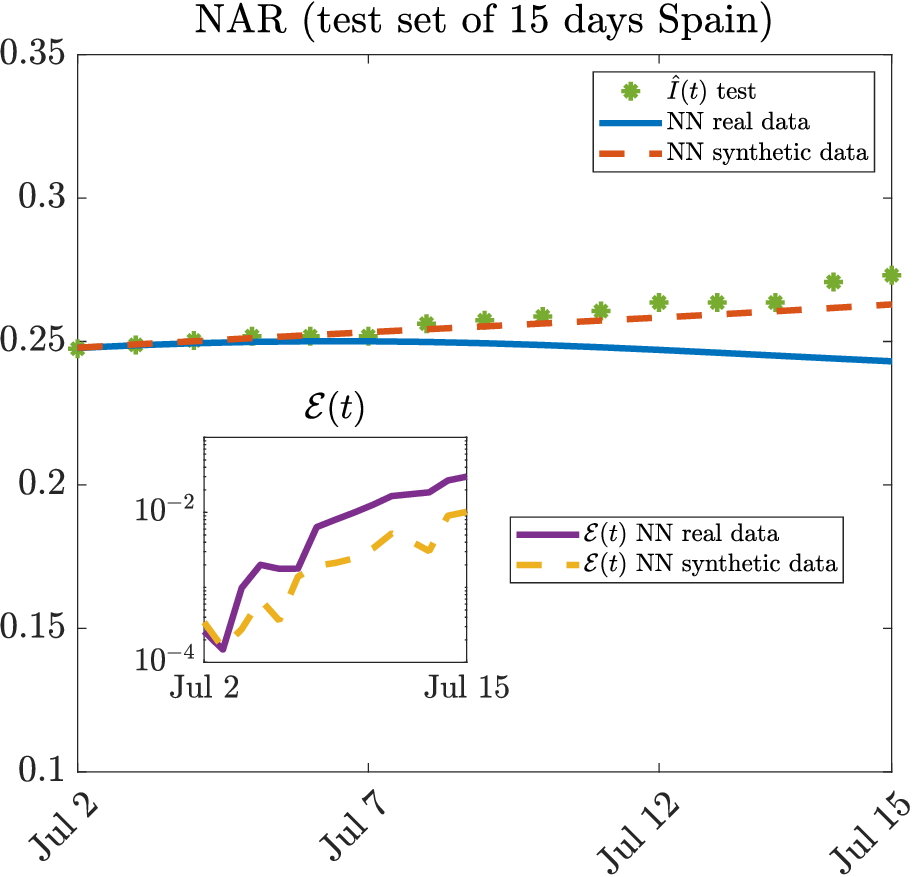}
	\caption{Nonlinear autoregressive network: short and long time predictions. Solution obtained by training a NAR network on both real and synthetic data, then testing it to produce daily data from July 2nd to July 8th (first row), and extend to July 15th (second row), compared with the available true data.     On the left, the case for Italy, and on the right, the one for Spain.}
	\label{fig:NAR_test}
\end{figure}

\paragraph{Nonlinear autoregressive neural network.}
We design now a NAR network with three fully connected layers, each consisting of $10$ nodes, and use the \texttt{tanh} activation function. To implement the neural network we use the Matlab function \texttt{narnet}. We consider a time interval $[t_l,T]$, being $t_l$ the time in which the lockdown started and $T$ the final time supposed to be July 1st. We discretize the interval with time step $\Delta t_r = 1$ if we deal with real data, or with time step $\Delta t_s = 0.2$ if we deal with synthetic data as before. The input to the network is the number of infected at time $t-2\Delta t,t- \Delta t$, and the output is $I(t)$ that is the number of infected at time $t$. Both the input and output data are generated as described above by considering a social-SIR model. As for the FNN we first evaluate the performance of the NAR network on the training set. Since the number of infected individuals during the training interval is known, we can use an open-loop prediction approach using the true training set values to forecast the next time step's solution.
Figure \ref{fig:NAR_training} shows the solution obtained by training a NAR network on both real and synthetic data, then testing it on the same dataset and comparing it to the available data. On the left, we have the results for Italy, and on the right, for Spain. Both the NAR networks, trained on real and synthetic data, are able to accurately approximate the underlying dynamics when one considers the training interval.
We then test the NAR network for both short-term and long-term predictions. For short-term predictions, we compare the neural network solution with a test set containing daily data from July 2nd to July 8th for both Italy and Spain, while for long-term predictions, we extend the test period to July 15th. Since these data are assumed to be unknown a priori, we convert the open-loop NAR network into a closed-loop form. This allows the network to generate predictions for the desired number of steps by using its previous predictions as input, without relying on the true values for evolution. Figure \ref{fig:NAR_test} shows the solution produced by the NAR network trained on both real and synthetic data, compared with the available true data. Additionally, we plot the pointwise error in time computed as in \eqref{eq:error}. In the first row, the short-term prediction while in the second row, the long term prediction.  On the left, Italy. On the right, Spain. The NAR network trained on synthetic data is able to provide a very accurate approximation of the solution and demonstrates strong predictive capabilities. Notably, it significantly outperforms the results obtained when training exclusively on real data, as it is particularly evident in the case of Italy.

\section{Conclusion}\label{sec:conclusion}
In this study, we explored the impact of a data augmentation strategy on enhancing the forecasting capabilities of Feed-Forward Neural Networks  and Nonlinear Autoregressive Networks in analyzing and predicting the dynamics of the COVID-19 epidemic during the lockdown and post-lockdown phase in Italy and Spain. By integrating real data with synthetic data generated from a generalized SIR model, we demonstrated the effectiveness of neural networks in capturing complex epidemic dynamics. FNNs performed well in interpolating within observed data, while NARs were particularly effective in capturing temporal dependencies, making them suitable for long-term forecasting. The generalized SIR model provided a robust foundation for generating high-resolution synthetic data, enriching the neural network training process. The key outcome of this work is the development of a hybrid methodology that bridges the gap between traditional mathematical modeling and machine learning techniques, such as Physics-Informed Neural Networks, offering a scalable and flexible framework for epidemic modeling and prediction in real-world scenarios.

While this study has focused on FNNs and NARs, other neural network architectures could be explored to further improve epidemic forecasting. For example, Long Short-Term Memory (LSTM) networks and Transformer-based models have been increasingly used for time-series prediction and could provide additional advantages in learning long-term dependencies \cite{amendolara2023lstm, nabi2021comparative, bontempi2023machine}. Similarly, hybrid AI-mechanistic models, which integrate deep learning with traditional compartmental models, have shown promise in epidemic modeling \cite{wang2021hybrid}. Future research could investigate whether these approaches provide added value in capturing epidemic trends compared to the models used in this study, particularly in scenarios with limited or noisy data.

An important direction for future research is the incorporation of uncertainty quantification (UQ) to better account for variability and noise in real-world data^^>\cite{bertaglia2022bi_fidelity, albi2022kinetic, capaldi2012parameter, chowell2017fitting}. While neural networks provide powerful forecasting tools, their deterministic nature limits their ability to express confidence in predictions, which is crucial in public health decision-making. Approaches such as Bayesian deep learning, Monte Carlo Dropout, and ensemble learning can be explored to quantify the range of possible epidemic trajectories. Additionally, integrating probabilistic surrogate models with mechanistic epidemiological frameworks can improve risk assessment and early-warning systems for future outbreaks. By embedding UQ techniques within the proposed hybrid framework, epidemic forecasting can become more robust, interpretable, and reliable to better account for variability and noise in real-world data.

\section*{Acknowledgements}
This work has been written within the activities of GNCS and GNFM groups of INdAM (Italian National Institute of High Mathematics). G.D. has been partially funded by the European Union — NextGenerationEU, MUR–PRIN 2022 PNRR Project No. P2022JC95T “Data-driven discovery and control of multi-scale interacting artiﬁcial agent systems”. G.D. and F.F. thank the Italian Ministry of University and Research (MUR) through the PRIN
2020 project (No. 2020JLWP23) ``Integrated Mathematical Approaches to Socio–Epidemiological Dynamics”. L.P. has been partially funded by the European Union– NextGenerationEU under the program “Future Artificial Intelligence– FAIR” (code PE0000013), MUR PNRR, Project “Advanced MATHematical methods for Artificial Intelligence– MATH4AI”.  L.P. acknowledges the support by the Royal Society under the Wolfson Fellowship ``Uncertainty quantification, data-driven simulations and learning of multiscale complex systems governed by PDEs" and by MIUR-PRIN 2022 Project (No. 2022KKJP4X), ``Advanced numerical methods for time dependent parametric partial differential equations with applications". The partial support by ICSC -- Centro Nazionale di Ricerca in High Performance Computing, Big Data and Quantum Computing, funded by European Union -- NextGenerationEU is also acknowledged. 
 
	\bibliographystyle{abbrv}
\bibliography{biblio}
\end{document}